\documentclass{article}
\usepackage[utf8]{inputenc}
\usepackage{amsmath}
\usepackage{csquotes}
\usepackage{amsthm}
\usepackage[
]{biblatex}
\usepackage{lipsum}
\usepackage{enumitem}
\usepackage{amsfonts}
\usepackage[english]{babel}
\usepackage{amsopn}

\newtheorem{theorem}{Theorem}
\newtheorem{corollary}{Corollary}[theorem]
\newtheorem{definition}{Definition}
\newtheorem{proposition}{Proposition}
\newtheorem{lemma}{Lemma}
\newtheorem{remark}{Remark}

\newtheorem{example}{Example}

\newcommand\blfootnote[1]{%
  \begingroup
  \renewcommand\thefootnote{}\footnote{#1}%
  \addtocounter{footnote}{-1}%
  \endgroup
}

\title{Long-term behavior of curve shortening flow \\in $\mathbb{R}^3$}

\author{Ji{\v r}{\' i} Minar{\v c}{\' i}k \and Michal Bene{\v s}}

\date{March 16, 2020
\footnote{Published in SIAM Journal on Mathematical Analysis, Vol. 52, Iss. 2 (2020).}
}

\begin{document}

\maketitle

\begin{abstract}
Space curve motion describes dynamics of material defects or interfaces, can be found in image processing or vortex dynamics. This article analyses some properties of space curves evolved by the curve shortening flow. In contrast to the classical case of shrinking planar curves, space curves do not obey the Avoidance principle in general. They can lose their convexity or develop non-circular singularities even if they are simple. In the first part of the text, we show that even though the convexity of space curves is not preserved during the motion, their orthogonal projections remain convex. In the second part, the Avoidance principle for spherical curves under the curve shortening flow in $\mathbb{R}^3$ is shown by generalizing the arguments developed by Hamilton and Gage.
\end{abstract}

\blfootnote{\url{https://doi.org/10.1137/19M1248522}}
\blfootnote{This work was supported by the Ministry of Education, Youth and Sports of the Czech Republic under the OP RDE grant number CZ.02.1.01/0.0 /0.0/16\_019/0000753 "Research centre for low-carbon energy technologies".}
\blfootnote{Department of Mathematics, Faculty of Nuclear Sciences and Physical Engineering, Czech Technical University in Prague, Prague 12000, Czech Republic.}


\section{Introduction}

Motion of curves and hypersurfaces induced by their mean curvature has been extensively studied both for its convenient mathematical properties and for its usefulness in applications ranging from physics to computer science. Adding an external forcing term to the prescribed normal velocity or considering an anisotropic geometrical setting leads to motion laws that can describe the dynamics of physical interfaces between different phases of material, defects in their crystalline structure or boundaries of thin layers (see [26, 18]).
Alternatively, such motion may be employed in the image segmentation task, where edges of different objects within a given image are automatically extracted and used for further processing (see [16, 6, 25]).
Similar motion laws for curves in $\mathbb{R}^3$ can describe motion of elastic rods discussed in \cite{Wardetzky2008}, evolution of vortex filaments via the localized induction approximation (see \cite{Ricca1991}) or dynamics of particle accumulation systems (see e.g. \cite{Schwabe2007}).

The theoretical efforts to understand properties of the original curve shortening problem in $\mathbb{R}^2$ have led to several important result obtained by Hamilton, Gage and Grayson in \cite{GaHa1986, Gray1987}. The well-known Grayson Theorem states that the curve shortening flow shrinks all simple planar curves to a point, making them asymptotically circular as they approach the singularity and keeping them simple throughout the timespan of the evolution. The last property is referred to as the Avoidance principle for planar curves (see \cite{White2002}). This property is particularly interesting as the normal velocity only depends on local geometrical information at a given point along the curve. The motion also preserves convexity of the curve and makes initially non-convex curves convex in finite time.

Many of classical results have been generalized for the mean curvature flow of hypersurfaces (see e.g. \cite{Huisken84}). However, these results do not hold for the codimension-two problems discussed in this contribution. Evolving curves embedded in the three-dimensional Euclidean space may develop local singularities before shrinking to a point and, in general, do not obey the Avoidance principle. Literature concerning the curvature flow of manifolds with higher codimension is rather sparse (see e.g. \cite{Zoo2013,Alt1991,AltGr1992,Corr16,Sim12,Khan15,MinKimBen18}). 

This problem was first studied by Altschuler and Grayson in \cite{Alt1991, AltGr1992}, where the short-term existence and uniqueness of the solution was shown. The article \cite{Alt1991} also classified all types of singularities that may develop during the motion. Recently, properties of this flow were studied in \cite{MinKimBen18, Khan15, Corr16, Sim12} and solitons of the flow were discussed in \cite{Zoo2013}.

This paper further addresses the discrepancies between the classical curve shortening flow in plane and the generalized codimension-two flow and aims to contribute to the understanding of long-term behavior of this problem. Convexity of curves and their two-dimensional projections during the flow is discussed in the first part of the text. The second part deals with spherical curves. First, we show that they obey the Avoidance principle and then discuss behavior of several spherical curves evolving at once.

The article is organized as follows. Section \ref{sec:codimension-two_curve_shortening_flow} introduces the necessary notation for the parametric formulation of the curve shortening problem in $\mathbb{R}^3$. The main results of this article are Proposition \ref{thm:projectable} and Proposition \ref{thm:spherical} in Section \ref{sec:convex_space_curves} and Section \ref{sec:spherical_curves}, respectively. In Section \ref{sec:convex_space_curves}, we consider a generalization of convexity for curves embedded in spaces of higher dimension and investigate whether it is preserved during the evolution. Then, the convexity of orthogonal projection into a given plane is analyzed in Proposition \ref{thm:projectable}. Section \ref{sec:spherical_curves} is focused on the evolution of spherical curves. It includes the proof of the Avoidance principle for spherical curves (see Proposition \ref{thm:spherical}). The last section contains final remarks and discussion of open problems.

\section{Curve shortening flow in $\mathbb{R}^3$} \label{sec:codimension-two_curve_shortening_flow}

This section introduces the curve shortening problem in $\mathbb{R}^3$ and necessary notation used in the article.

The problem can be formalized in several ways. In this contribution, we use the parametric approach for its simplicity. Other approaches, such as the phase field or the level set method (see \cite{Seth81, SethOsher88}), are better suited for planar curves. Generalization of these methods for manifolds with higher codimension can be found in \cite{BurchChengMerrOsher2001} or \cite{AmbrosSoner96}.

Let $\{ \Gamma_t \}_{t \in [0, t_{max})}$ be a family of closed curves in $\mathbb{R}^3$ evolving in time interval $[0, t_{max})$. Each curve $\Gamma_t$ is given by a parametric function $X(\cdot,t) : \mathbb{S}^1 \rightarrow \mathbb{R}^3$, where $\mathbb{S}^1=\mathbb{R}/2\pi\mathbb{Z}$ is a unit circle. We assume that the rate of parametrization does not vanish, i.e. $\Vert \partial_u X(u,t) \Vert > 0$ for all $u \in \mathbb{S}^1$ and $t \in [0,t_{max})$, where $\Vert \cdot \Vert$ is the Euclidean norm in $\mathbb{R}^3$. Furthermore, we require $X$ to be a $\mathcal{C}^1$-class function and $X(\cdot, t)$ to be a $\mathcal{C}^2$-class function for any given $t \in [0,t_{max})$.

The local geometry of $\Gamma_t$ at each point $u \in \mathbb{S}^1$ is given by the Frenet frame $\{ T(u,t), N(u,t), B(u,t) \}$, the curvature $\kappa(u,t)$ and the torsion $\tau(u,t)$. Using the symbol $s$ to denote the arclength parametrization which satisfies $\Vert \partial_s X \Vert = 1$, we can define the tangent vector $T$ and the curvature $\kappa$ as $T = \partial_s X$ and $\kappa = \Vert \partial_s^2 X \Vert$, respectively. The normal vector $N$ exists only when $\kappa \neq 0$ and is given by $N = \frac{1}{\kappa} \partial_s^2 X$. The binormal vector $B$ is then defined as the cross product between $T$ and $N$.

The curve shortening flow is defined as the following initial-value problem for the parametrization $X = X(u,t)$:
\begin{align}
\partial_t X &= \kappa N & \mbox{on } \mathbb{S}^1 \times [0,t_{max}), \label{eq:curvature_flow_1} \\
X\vert_{t = 0} &= X_0 & \mbox{in } \mathbb{S}^1, \label{eq:curvature_flow_2}
\end{align}
where $X_0 \in \mathcal{C}^2 (\mathbb{S}^1; \mathbb{R}^3)$ is the parametrization of the initial curve $\Gamma_0$. Although $N$ is undefined for all points on $\Gamma_t$ where $\kappa = 0$, the term $\kappa N = \partial_s^2 X$ in (\ref{eq:curvature_flow_1}) remains defined everywhere.

\section{Convex space curves} \label{sec:convex_space_curves}

As mentioned in the introduction, the curve shortening flow preserves convexity of planar curves. The aim of this section is to find an analogous statement for the codimension-two motion (\ref{eq:curvature_flow_1}-\ref{eq:curvature_flow_2}).

The notion of convexity for manifolds of higher codimension is not commonly defined. In \cite{Sedykh1994}, space curves are called convex if they lie on the boundary of their convex hull. This definition becomes troublesome for planar curves where the boundary of the convex hull is equal to the hull itself. To address this issue, we propose the following definition employing the Minkowski functional.

\begin{definition}[Convex space curve]
\label{def:convex}
For any convex set $K \subset \mathbb{R}^3$, let $M_K : \mathbb{R}^n \rightarrow \mathbb{R}^+$ denote the Minkowski functional prescribed by
\begin{align*}
M_K(x) := \inf \left\{ \lambda \in \mathbb{R}^+ \; : \; \tfrac{1}{\lambda} x \in K \right\}.
\end{align*}
for all $x \in \mathbb{R}^3$. We say that a closed curve $\Gamma$ is convex if $M_{\mathrm{C}(\Gamma)}|_{\Gamma} \equiv 1$, where $\mathrm{C}(\Gamma)$ is the convex hull of $\Gamma$.
\end{definition}

The following example illustrates that in contrast to the planar case, convex space curves may loose their convexity during the evolution by (\ref{eq:curvature_flow_1}-\ref{eq:curvature_flow_2}).

\begin{example}\label{ex:counter-example-convex}
Consider a convex curve $\Gamma_0 = \mbox{Ran} \, X_0$ given by
\begin{align}
\label{eq:counter-example-convex}
X_0(u) := 
\begin{pmatrix} 
\cos ( au^3 + bu) \\
\sin ( au^3 + bu) \\
\sin u - \tfrac{1}{2} \sin(2u)
\end{pmatrix}
\end{align}
for $u \in \mathbb{S}^1$, where $a = \frac{\pi + 2}{2(1-\pi^2)}$ and $b = \frac{\pi^3 + 2}{2(\pi^2 - 1)}$. Let $\{ \Gamma_t \}_{t \in[0,t_{max})}$ be a family of space curves evolving according to the curve shortening flow given by (\ref{eq:curvature_flow_1}-\ref{eq:curvature_flow_2}) with the initial condition $\Gamma_0$. Since the curvature of the original curve $\Gamma_0$ at the point $u = 0$ is greater than its curvature at $u = w$ and $u = -w$, where $w = (\pi^3 + 2)^{\frac{1}{2}}(\pi + 2)^{-\frac{1}{2}}$, $X(0,t)$ departs from the line segment $\mathrm{C}(\{ X(-w,t), X(w,t) \})$ which lies on the boundary of $\mathrm{C} (\Gamma_t)$. Thus $\Gamma_t$ will stop being convex immediately after $t = 0$.
\end{example}

Example \ref{ex:counter-example-convex} shows that the convexity proposed in Definition \ref{def:convex} is not preserved. In Proposition \ref{thm:projectable}, we show that even though space curves may lose their convexity, their initially convex orthogonal projections remain convex throughout the total timespan of the evolution.

The following two lemmas are used in the proof of Proposition \ref{thm:projectable}. Lemma \ref{lem:star-shaped} uses the notion of star-shaped curves, which are boundaries of star-shaped sets. In Lemma \ref{lem:star-shaped}, we use the following sufficient condition. If $\Gamma$ is a closed planar curve and there is $x \in \mathrm{int} \, \Gamma$ such that $X(u,t) - x$ and $T(u,t)$ are linearly independent for all $u \in \mathbb{S}^1$, then $\Gamma$ is star-shaped.

\begin{lemma}\label{lem:star-shaped}
Let $\{ \Gamma_t \}_{t \in [0,t_{max})}$ be a family of closed planar curves such that $\Gamma_0$ is convex and the parametrization $X \in \mathcal{C}^1(\mathbb{S}^1 \times [0, t_{max}))$. Then there exists $t_0 \in (0,t_{max})$ such that $\Gamma_t$ is star-shaped for all $t \in [0,t_0)$.
\begin{proof}
Select any $x$ from $\mathrm{int} \, \Gamma_0$ and define a function
\begin{align*}
\varphi (u,t) := \Vert X(u,t) - x \Vert - \vert \langle X(u,t) - x, T(u,t) \rangle \vert,
\end{align*}
where $\varphi \in \mathcal{C}(\mathbb{S}^1 \times [0, t_{max}))$ due to the assumptions. Furthermore, $\varphi$ is non-negative because of the Cauchy-Schwarz inequality. Next, $\inf_{\mathbb{S}^1} \varphi|_{t = 0} > 0$ because $\Gamma_0$ is convex and thus star-shaped with respect to any inner point. Since $\varphi$ is continuous, there is $t_0 \in [0,t_{max})$ which satisfies
\begin{align*}
\forall t \in [0,t_0): \; x \in \mathrm{int}\, \Gamma_t \;\;\; \wedge \;\;\; \inf_{\mathbb{S}^1\times[0,t_0)} \varphi > \tfrac{1}{2} \inf_{\mathbb{S}^1} \varphi|_{t = 0}. 
\end{align*}
This implies the original statement.
\end{proof}
\end{lemma}

Notice that Lemma \ref{lem:star-shaped} imposes regularity assumption on the parametrization $X$ but it does not require the curve to follow the curve shortening equation (\ref{eq:curvature_flow_1}-\ref{eq:curvature_flow_2}).

\begin{lemma}\label{lem:PN-N_P}
Let $\mathnormal{P} \in \mathcal{L}(\mathbb{R}^3)$ be an orthogonal projection with $\mathnormal{\emph{dim}} \mathnormal{\emph{Ran}} \, \mathnormal{P} = 2$ and $\Gamma$ be a space curve in $\mathbb{R}^3$. Assume that $\Vert \partial_u \mathnormal{P} \mathnormal{X}(u) \Vert > 0$ and $\kappa_P(u) > 0$ for some $u \in \mathbb{S}^1$, where $\kappa_P$ denotes the curvature of the projected curve $\mathnormal{P} \Gamma$ at point $PX(u)$. Then
\begin{align*}
\left\langle \mathnormal{P} \mathnormal{N}(u) , \mathnormal{N}_{\mathnormal{P}}(u) \right\rangle > 0,
\end{align*}
where $\mathnormal{N}(u)$ is the normal vector of $\Gamma$ at the point $\mathnormal{X}(u)$ and $\mathnormal{N}_{\mathnormal{P}}(u)$ denotes the normal vector of the projected curve $\mathnormal{P}\Gamma$ at point $\mathnormal{P} \mathnormal{X}(u)$.
\begin{proof}
The assumptions imply that $\kappa(u) > 0$ and that both $N$ and $N_P$ are well defined at $u$. Note that we omit explicitly writing the argument $u$ in the rest of the proof to increase its readability. Using $\partial_s = \Vert \partial_u X \Vert^{-1} \partial_u$ to denote the arclength derivative with respect to the original curve $\Gamma$, we obtain
\begin{align*}
\left\langle \mathnormal{P} \mathnormal{N}, \mathnormal{N}_{\mathnormal{P}} \right\rangle &= \left\langle \frac{1}{\kappa} \partial_s^2 \mathnormal{P} \mathnormal{X}, \frac{1}{\kappa_{\mathnormal{P}} \Vert \partial_s \mathnormal{P} \mathnormal{X} \Vert^2} \partial_s^2 \mathnormal{P} \mathnormal{X} - \frac{\partial_s \Vert \partial_s \mathnormal{P} \mathnormal{X} \Vert}{\kappa_{\mathnormal{P}} \Vert \partial_s \mathnormal{P} \mathnormal{X} \Vert^3} \partial_s \mathnormal{P} \mathnormal{X} \right\rangle \\
&= \frac{\Vert \partial_s^2 \mathnormal{P} \mathnormal{X} \Vert^2 \Vert \partial_s \mathnormal{P} \mathnormal{X} \Vert^2 - \langle \partial_s \mathnormal{P} \mathnormal{X}, \partial_s^2 \mathnormal{P} \mathnormal{X} \rangle^2}{\kappa \kappa_{\mathnormal{P}} \Vert \partial_s \mathnormal{P} \mathnormal{X} \Vert^4}.
\end{align*}
The inequality $\left\langle PN, N_P \right\rangle \geq 0$ is obtained directly from the Cauchy-Schwarz inequality
\begin{align}\label{ineq:new-Cauchy-Schwarz}
\langle \partial_s \mathnormal{P} \mathnormal{X}, \partial_s^2 \mathnormal{P} \mathnormal{X} \rangle \leq \Vert \partial_s \mathnormal{P} \mathnormal{X} \Vert \Vert \partial_s^2 \mathnormal{P} \mathnormal{X} \Vert.
\end{align}
Equality can occur only when there is some $\alpha \geq 0$ such that $\partial_s^2 PX = \alpha \partial_s PX$. Note that $\partial_s PX$ is non-zero from the assumptions. This implies that $\partial_u^2 PX = \beta \partial_u PX$, where $\beta = \alpha \Vert \partial_u X \Vert + \tfrac{1}{2} \partial_u \Vert \partial_u X \Vert^2$. Thus we obtain
\begin{align*}
\kappa_P N_P &= \Vert \partial_u PX \Vert^{-2} \partial_u^2 PX - \Vert \partial_u PX \Vert^{-3} \partial_u \Vert \partial_u PX \Vert \partial_u PX \\
&= \beta \Vert \partial_u PX \Vert^{-2} - \Vert \partial_u PX \Vert^{-4} \left\langle \partial_u PX, \partial_u^2 PX \right\rangle = 0.
\end{align*}
Since $\kappa_P$ is assumed to be positive, $\partial_s PX$ and $\partial_s^2 PX$ must be linearly independent and the inequality \ref{ineq:new-Cauchy-Schwarz} is strict.
\end{proof}
\end{lemma}

As shown in \ref{ex:counter-example-convex}, space curves may lose their convexity during the motion. The following proposition states that the convexity of their orthogonal projections is preserved.

\begin{proposition}\label{thm:projectable}
Let $\mathnormal{P} \in \mathcal{L}(\mathbb{R}^3)$ be an orthogonal projection with $\mathnormal{\emph{dim}} \, \mathnormal{\emph{Ran}} \, \mathnormal{P} = 2$ and $\Gamma_0$ be a space curve such that its projection $P\Gamma_0$ is convex. Assume that the parametrization $PX(\cdot, t)$ of the projected curve $P\Gamma_t$ is regular for all $t \in [0, t_{max})$, i.e. $PT$ does not vanish on $\mathbb{S}^1 \times [0, t_{max})$. If $\Gamma_t$ evolves according to the curve shortening flow given by \ref{eq:curvature_flow_1}\ref{eq:curvature_flow_2} with the initial condition $\Gamma_0$, then $P\Gamma_t$ is convex for all $t \in [0,t_{max})$.
\begin{proof}
Assume that $P\Gamma_t$ looses its convexity during the evolution. In order to formalize the proof, we define the following auxiliary functional:
\begin{align*}
\phi &: \mathbb{S}^1\times[0,t_{max}) \rightarrow \mathbb{R}^+_0: \; (u,t) \mapsto \mbox{dist} \left( PX(u,t), \;\; \partial \mathrm{C}( P\Gamma_t ) \right).
\end{align*}
The proof is divided into several steps in which individual statements (a), (b), (c) and (d) are shown. Their combination then leads to a contradiction.

\begin{enumerate}[label=(\alph*)]
\item{$\forall t \in [0,t_{max}): \;\; \phi(\cdot, t)$ is continuous on $\mathbb{S}^1$.}
\end{enumerate}

\noindent We know that $X(\cdot,t)$ is continuous, i.e. for all $t \in [0,t_{max})$, $u \in \mathbb{S}^1$ and $\varepsilon > 0$ there is $\delta > 0$ such that $\vert u' - u \vert < \delta$ implies $\Vert X(u', t) - X(u,t) \Vert < \varepsilon$ for all $u' \in \mathbb{S}^1$. This allows us to write
\begin{align*}
\phi(u',t) &= \inf\limits_{Y \in \partial \mathrm{C}( P\Gamma_t ) } \Vert PX(u',t) - Y \Vert \\
&\leq \Vert PX(u',t) - PX(u,t) \Vert + \inf\limits_{Y \in \partial \mathrm{C}( P\Gamma_t ) } \Vert PX(u,t) - Y \Vert \\
&\leq \Vert X(u',t) - X(u,t) \Vert + \phi(u,t) < \varepsilon + \phi(u,t).
\end{align*}
Similarly $\phi(u,t) < \varepsilon + \phi(u',t)$ and thus $|\phi(u',t) - \phi(u,t)| < \varepsilon$ when $|u' - u| < \delta$. 

Because $\mathbb{S}^1$ is compact and $\phi(\cdot, t)$ is continuous on $\mathbb{S}^1$, $\phi(\cdot,t)$ attains its maximum on $\mathbb{S}^1$. The maximum at each time $t$ is denoted by an auxiliary function $\Phi$.
\begin{align*}
\Phi &: [0,t_{max}) \rightarrow \mathbb{R}^+_0: \; t \mapsto \max\limits_{u \in \mathbb{S}^1} \phi(u,t).
\end{align*}
Note that $P\Gamma_t$ is convex if and only if $\Phi(t) = 0$.

\begin{enumerate}[label=(\alph*)]
\setcounter{enumi}{1}
\item{$\Phi$ is continuous on $[0,t_{max})$.}
\end{enumerate}

\noindent It suffices to show that $\phi$ is continuous with respect to time $t$. This would mean that for all $t \in [0,t_{max})$, $u \in \mathbb{S}^1$ and $\varepsilon > 0$, there is $\delta > 0$ such that for all $t' \in [0,t_{max})$ and all $u'$, $\vert u' - u \vert < \delta$ implies $\vert \phi(u,t') - \phi(u,t) \vert < \varepsilon$. Since $\partial \mathrm{C}( P\Gamma_t )$ is compact, there exists $\tilde{Y} \in \partial \mathrm{C}( P\Gamma_t )$ such that
\begin{align}
\phi(u,t) = \inf\limits_{Y \in \partial \mathrm{C}( P\Gamma_t )} \Vert PX(u,t) - Y \Vert = \Vert PX(u,t) - \tilde{Y} \Vert.
\end{align}
For $t'$ close enough to $t$, the set $\partial \mathrm{C}( P\Gamma_t )$ is similar to $\partial \mathrm{C}( P\Gamma_{t'} )$ in terms of the Hausdorff distance. This means that there is $\tilde{Y}' \in \partial \mathrm{C}( P\Gamma_{t'} )$ such that $\Vert \tilde{Y}' - \tilde{Y} \Vert$ can be arbitrarily small if $t'$ and $t$ are close enough. Then
\begin{align*}
\phi(u,t') &= \inf\limits_{Y \in \partial \mathrm{C}( P\Gamma_{t'} ) } \Vert PX(u,t') - Y \Vert \leq \Vert PX(u,t') - \tilde{Y}' \Vert \\
&\leq \Vert PX(u,t') - PX(u,t) \Vert + \Vert PX(u,t) - \tilde{Y} \Vert + \Vert \tilde{Y} - \tilde{Y}' \Vert \\
&= \Vert X(u,t') - X(u,t) \Vert + \Vert \tilde{Y} - \tilde{Y}' \Vert + \phi(u,t) < \varepsilon + \phi(u,t).
\end{align*}
Similarly $\phi(u,t) < \varepsilon + \phi(u,t')$. Therefore $|\phi(u,t') - \phi(u,t)| < \varepsilon$ when $|t' - t| < \delta$.

Let $t_0$ denote the time when $\Gamma$ looses its convexity, i.e.
\begin{align}
\label{eq:t_0}
t_0 := \inf\{ t \in (0,t_{max}) : \; \Phi(t) > 0 \},
\end{align}
where the set $\{ t \in (0,t_{max}) : \; \Phi(t) > 0 \}$ is nonempty by the assumption.

\begin{enumerate}[label=(\alph*)]
\setcounter{enumi}{2}
\item{$\Phi(t_0) = 0$.}
\end{enumerate}

\noindent If $\Phi(t_0) > 0$, continuity of $\Phi$ on $[0,t_{max})$ implies
\begin{align*}
\exists \varepsilon > 0, \;\; \forall t \in (t_0 - \varepsilon, t_0): \;\; \Phi(t) > \tfrac{1}{2}\Phi(t_0)>0.
\end{align*}
This would contradict the definition of $t_0$.

\begin{enumerate}[label=(\alph*)]
\setcounter{enumi}{3}
\item{$\exists t_1 \in (t_0,t_{max}), \;\; \forall t \in (t_0, t_1): \;\; \Phi(t)$ is non-increasing.}
\end{enumerate}

\noindent By the assumption, the projection $P\Gamma_t$ is a regular curve. Since $P$ is a linear operator, the regularity of the parametrization $PX$ is at least $\mathcal{C}^1$. Thus, we may use \ref{lem:star-shaped}, which states that there exists $t_1 \in (t_0,t_{max})$ such that $P\Gamma_t$ is star-shaped for all $t \in (t_0,t_1)$. For time $t \in (t_0,t_1)$ fixed, the function $\phi(\cdot,t)$ reaches its maximum at the point denoted by $u_2 \in \mathbb{S}^1$. Let $u_1, u_3 \in \mathbb{S}^1$ satisfy $u_2 \in [u_1, u_3]$, $PX(u_1, t) \in \partial \mathrm{C} (P \Gamma_t)$, $PX(u_3, t) \in \partial \mathrm{C} (P \Gamma_t)$ and
\begin{align*}
\forall u \in (u_1, u_3) : \; PX(u,t) \notin \partial \mathrm{C} (P \Gamma_t).
\end{align*}

Consider orthogonal coordinate system $x$-$y$ in $\mbox{Ran} \, P$ such that the $x$-axis is parallel to $PX(u_3,t) - PX(u_1,t)$ and $PX(u_2,t) - PX(u_1,t)$ has a positive $y$ coordinate. Since $P \Gamma_t$ is star-shaped, no small kinks can develop along the curve and we can express $\Phi(t)$ as the difference between the $y$ coordinate of $PX(u_2,t)$ and $PX(u_1,t)$.

If $\kappa_P(u_2, t) > 0$, \ref{lem:PN-N_P} and the motion law \ref{eq:curvature_flow_1} imply that the $y$ coordinate of $PX(u_2, t)$ is non-increasing in time. Similarly, when $\kappa_P(u_1, t) > 0$ and/or $\kappa_P(u_3, t) > 0$, the $y$ coordinate of $PX(u_1, t)$ and/or $PX(u_3, t)$ is non-decreasing, respectively. When $\kappa_P = 0$ at $u_1$, $u_2$ or $u_3$, the $y$ coordinate of the corresponding point remains constant as the motion takes place only in the $x$ direction. In any case, the distance between $PX(u_2, t)$ and $\partial \mathrm{C}(P\Gamma_t)$ cannot increase and thus $\Phi(t)$ cannot increase either.

Finally, using (b), (c) and (d), we get
\begin{align*}
\forall t \in (t_0,t_1): \;\; 0 \leq \Phi(t) \leq \Phi(t_0) = 0,
\end{align*}
which implies that $\Phi = 0$ in $(t_0,t_1)$. This contradicts the definition of $t_0$ in \ref{eq:t_0} and thus ensures the convexity of $P\Gamma_t$ for all $t \in [0,t_{max})$.
\end{proof}
\end{proposition}

\section{Spherical curves} \label{sec:spherical_curves}
The space curve $\Gamma$ is called spherical if there exists a point $x \in \mathbb{R}^3$ and a positive constant $r$ such that $\Vert x - y \Vert = r$ for all $y \in \Gamma$. Thanks to their convenient properties, spherical curves and their behavior during the shortening flow has gained new attention in recent years. It has been recently discovered in \cite{Sim12} that initially spherical curves remain spherical during the flow. We refer the reader to \cite{Khan15} for further discussion and consequences of this result and to \cite{MinKimBen18} for an alternative proof.

In this section, we show that spherical curves also satisfy the Avoidance principle and thus closely resemble the behavior of evolving planar curves. This result is achieved by generalizing the classical proof from \cite{GaHa1986}. Similar process might be used for other types of space curves. Main obstacle is solved by the following lemma.

\begin{lemma}\label{lem:collinear}
Let $\Gamma$ be a space curve embedded in a sphere and denote $f(u_1, u_2) := \Vert X(u_2,t) - X(u_1,t) \Vert^2$ for all $u_1, u_2 \in \mathbb{S}^1$. If $f$ has a local minimum at $( u_1 ,  u_2)$ and $ u_1  \neq  u_2$, then $T( u_1 )$ and $T( u_2)$ are collinear. 
\begin{proof}
Since $f$ has an extremum at $( u_1 , u_2)$, we get
\begin{align*}
\nabla f(u_1 , u_2) &= 2
\begin{pmatrix}
\Vert \partial_u X(u_1) \Vert \langle X(u_1) - X(u_2), T(u_1) \rangle \\
\Vert \partial_u X(u_2) \Vert \langle X(u_2) - X(u_1), T(u_2) \rangle
\end{pmatrix} = 
\begin{pmatrix}
0 \\ 0
\end{pmatrix}.
\end{align*}
Thus $T( u_1 )$ and $T( u_2 )$ are orthogonal to $X( u_2 ) - X( u_1 )$. Let $x$ denote the center of the sphere. Then $X(u_1) - x$ and $X(u_2) - x$ are orthogonal to the tangent plane of the sphere at the point $X(u_1)$ and $X( u_2)$, respectively. Thus $T(u_1)$ and $T(u_2)$ are orthogonal to $X(u_1) - x$ and $X(u_2) - x$, respectively. Together, we have
\begin{align*}
T( u_1 ) &\in \mathrm{span}(X( u_2 ) - X( u_1 ))^{\perp} \cap \mathrm{span}(X( u_1 ) - x)^{\perp}, \\
T( u_2 ) &\in \mathrm{span}(X( u_1 ) - X( u_2 ))^{\perp} \cap \mathrm{span}(X( u_2 ) - x)^{\perp},
\end{align*}
where $\mathrm{span} \, M$ and $M^{\perp}$ denote the linear span and the orthogonal complement of the set $M$, respectively. Then linearity of the inner product $\langle \cdot, \cdot \rangle$ implies that
\begin{align}\label{eq:intersection}
T( u_1 ), T( u_2 ) \in \mathrm{span}(X( u_1 ) - x)^{\perp} \cap \mathrm{span}(X( u_2 ) - x)^{\perp}.
\end{align}
Since $ u_1  \neq  u_2$, the intersection of $\mathrm{span}(X( u_1 ) - x)^{\perp}$ and $\mathrm{span}(X( u_2) - x)^{\perp}$ from \ref{eq:intersection} is a one-dimensional affine space. Thus $T( u_1 )$ and $T( u_2)$ are collinear.
\end{proof}
\end{lemma}

\begin{remark}
Note that \ref{lem:collinear} does not generalize to curves embedded in hyperspheres in $\mathbb{R}^n$ for $n > 3$. For example, consider the following curve in $\mathbb{R}^4$ with the parametrization
\begin{align*}
X(u) := \begin{pmatrix} \sin(\cos u) \\ \cos(\cos u) \sin (\sin(2u)) \\ \cos(\cos u)  \cos(\sin(2u)) \cos (\frac{1}{2} \sin u) \\ \cos(\cos u) \cos (\sin(2u)) \sin (\frac{1}{2} \sin u)
\end{pmatrix},
\end{align*}
and a tuple $(u_1,u_2) = \left(\frac{\pi}{2}, \frac{3\pi}{2}\right)$. This curve is embedded in unit 3-sphere, $f$ has a local minimum at the point $(u_1,u_2)$ and yet $\langle T(u_1), T(u_2) \rangle = 0$. 
\end{remark}

Following the proof from \cite{GaHa1986}, we define the functional
\begin{align}\label{eq:functional_f}
f : \mathbb{T} \times [0,t_{max}) \rightarrow \mathbb{R}^+: \; ( u_1 ,  u_2, t) \mapsto \Vert \mathnormal{X}( u_2, t) - \mathnormal{X}( u_1 , t) \Vert^2,
\end{align}
where $\mathbb{T} = \mathbb{S}^1 \times \mathbb{S}^1$ is a torus. The following lemma, inspired by \cite{GaHa1986}, trivially generalizes to our setting of spherical curves in $\mathbb{R}^3$.

\begin{lemma}[Lemma 3.2.2 from \cite{GaHa1986}] \label{lem:heat}
The functional $f$ defined in \ref{eq:functional_f} satisfies a strictly parabolic partial differential equation
\begin{align}
\label{eq:heat_eq_for_dp}
\partial_t f - \Delta f = -4,
\end{align}
where $\Delta = \partial^2_{s_1} + \partial^2_{s_2}$.
\begin{proof}
Using $\partial_t X = \kappa N$ and $\partial_s X = T$, a straightforward differentiation yields
\begin{align}
\partial_t f(u_1,u_2,t) &= 2 \langle X(u_1,t) - X(u_2,t) , \kappa(u_1,t)N(u_1,t) - \kappa(u_2,t)N(u_2,t) \rangle, \label{eq:dt_dp}\\
\partial_{s_1} f(u_1,u_2,t) &= 2 \langle X(u_1,t) - X(u_2,t) , T(u_1,t) \rangle, \nonumber \\
\partial_{s_2} f(u_1,u_2,t) &= -2 \langle X(u_1,t) - X(u_2,t) , T(u_2,t) \rangle, \nonumber \\
\partial^2_{s_1} f(u_1,u_2,t) &= 2 + 2 \langle X(u_1,t) - X(u_2,t) , \kappa(u_1,t) N(u_1,t) \rangle, \label{eq:ds1ds1_dp}\\
\partial^2_{s_2} f(u_1,u_2,t) &= 2 - 2 \langle X(u_1,t) - X(u_2,t), \kappa(u_2,t)N(u_2,t) \rangle. \label{eq:ds2ds2_dp}
\end{align}
Adding \ref{eq:dt_dp}, \ref{eq:ds1ds1_dp} and \ref{eq:ds2ds2_dp} proves \ref{eq:heat_eq_for_dp}.
\end{proof}
\end{lemma}

The following Lemma ensures the absence of small kinks in a curve with bounded curvature. The original result for planar curves is due to Schur \cite{Schur1921}. Later, the generalized version was given by Schmidt \cite{Schmidt1925}. The following formulation of Schur Theorem has been adopted from \cite{Lopez11}.

\begin{lemma}[Schur Comparison Theorem]\label{lem:schur}
Let $\Gamma_1$ and $\Gamma_2$ be open space curves of the same length $L \in \mathbb{R}^+$ with the arc-length parametrizations $X_1, X_2 : [0,L] \rightarrow \mathbb{R}^3$. Assume that $\Gamma_1$ is planar, and
\begin{align*}
\Gamma_1 \cup \partial \mathrm{C} ( \Gamma_1 ) \subset \partial \mathrm{C} (\Gamma_1 \cup \partial \mathrm{C} ( \Gamma_1 )).
\end{align*}
Furthermore, assume that $\kappa_1(s) \geq \kappa_2(s)$ for all $s \in [0,L]$, where $\kappa_1$ and $\kappa_2$ are curvatures of $\Gamma_1$ and $\Gamma_2$, respectively. Then $\Vert X_1(0) - X_1(L) \Vert \leq \Vert X_2(0) - X_2(L) \Vert$.
\end{lemma}

Using the generalized version of Schur Comparison Theorem, stated in \ref{lem:schur}, we can reformulate the following result from \cite{GaHa1986}.

\begin{lemma}[Corollary 3.2.4 from \cite{GaHa1986}]\label{lem:ineq}
Let $\Gamma$ be a space curve with uniformly bounded curvature, i.e.
\begin{align*}
\exists C > 0, \,\, \forall u \in \mathbb{S}^1, \,\, \forall t \in [0,t_{max}): \,\, \kappa(u, t) < C.
\end{align*}
Then the functional $f$ given by \ref{eq:functional_f} satisfies
\begin{align*}
f(u_1, u_2, t) \geq \frac{4}{C^2} \left[ \sin\left( \frac{2}{C} \int_{u_1}^{u_2} \Vert \partial_u X(u,t) \Vert \, \mathrm{d} u \ \right)\right]^2,
\end{align*}
for all $(u_1,u_2,t) \in \mathbb{T}^2 \times [0, t_{max})$.
\begin{proof}
The inequality is obtained from \ref{lem:schur} with $\Gamma_1$ being arc of a circle with radius $C^{-1}$. 
\end{proof}
\end{lemma}

The following proposition is the Avoidance principle for spherical curves. It states that originally simple spherical curves with bounded curvature cannot intersect themselves during the evolution. The statement and its proof are based on the classical result for planar curves from \cite{GaHa1986}.

\begin{proposition}[Avoidance principle for spherical curves]\label{thm:spherical}
Let $\Gamma_t$ evolve according to the curve shortening flow given by \ref{eq:curvature_flow_1}\ref{eq:curvature_flow_2} with the initial condition $\Gamma_0$. Assume that $\Gamma_0$ is a simple spherical curve and its curvature $\kappa$ can be uniformly bounded by a positive constant $C$ for all $u \in \mathbb{S}^1$ and $t \in [0,t_{max})$. Then $\Gamma_t$ cannot intersect itself for all $t \in (0,t_{max})$.
\begin{proof}
We split $\mathbb{T}^2 \times [0,t_{max})$ into $E$ and $D = \left[ \mathbb{T}^2 \times [0,t_{max}) \right] \setminus E$, where
\begin{align}\label{eq:def_E}
E := \left\{ (u_1,u_2,t) \in \mathbb{T}^2 \times [0,t_{max}) : \int_{u_1}^{u_2} \Vert \partial_u X(u,t) \Vert \, \mathrm{d} u < \frac{\pi}{C} \right\}.
\end{align}
\ref{lem:ineq} ensures the fact that $f(u_1,u_2,t) = 0$ implies $u_1 = u_2$ for all $(u_1,u_2,t) \in E$. In order to prove that $\Gamma_t$ cannot intersect itself, it suffices to show that
\begin{align*}
\inf_{D} f(u_1,u_2,t) > 0.
\end{align*}
From the definition of $E$ in \ref{eq:def_E}, all $(u_1,u_2,t)$ from the boundary $\partial E$ satisfy
\begin{align*}
\int_{u_1}^{u_2} \Vert \partial_u X(u,t) \Vert \, \mathrm{d} u = \frac{\pi}{C}.
\end{align*} 
Together with \ref{lem:ineq}, we have $f \geq \frac{4}{C^2}$ on $\partial E$. Since the initial curve $\Gamma_0$ is embedded and closed, there exists $m_1 \in (0, \frac{4}{C^2})$ such that
\begin{align*}
\inf_{\partial D}f(u_1,u_2,t) \geq \min \left\{ \inf_{\partial D \setminus \partial E}f(u_1,u_2,t), \frac{4}{C^2} \right\} > m_1.
\end{align*}
For $\varepsilon >0$ set
\begin{align}\label{eq:f_varepsilon}
f_\varepsilon(u_1,u_2,t) := f(u_1,u_2,t) + \varepsilon t.
\end{align}
Assume there exists $m_2 \in (0, m_1)$ and $(u_1^o,u_2^o,t^o) \in D$ such that $f_{\varepsilon} (u_1^o,u_2^o,t^o) = m_2$, where $t^o$ is the smallest possible. Since $f_{\varepsilon}$ attains its local minimum at $(u_1^o,u_2^o,t^o)$, we may use \ref{lem:collinear} and conclude that
\begin{align}\label{eq:rovna_se_2}
|\partial_{s_1} \partial_{s_2} f_{\varepsilon}(u_1^o,u_2^o,t^o)| = 2 | \langle T(u_1^o,t^o), T(u_2^o,t^o) \rangle |  = 2 \Vert T(u_1^o,t^o) \Vert \Vert T(u_2^o,t^o) \Vert = 2.
\end{align}
Since $t^o$ is the smallest possible, we get $\partial_t f_{\varepsilon}(u_1^o,u_2^o,t^o) \leq 0$ and
\begin{align} \label{eq:determinant}
\det \begin{pmatrix} \partial^2_{s_1} f_{\varepsilon}(u_1^o,u_2^o,t^o) & \partial_{s_1}\partial_{s_2} f_{\varepsilon}(u_1^o,u_2^o,t^o) \\ \partial_{s_2}\partial_{s_1} f_{\varepsilon}(u_1^o,u_2^o,t^o) & \partial^2_{s_2} f_{\varepsilon}(u_1^o,u_2^o,t^o) \end{pmatrix} &= \\
= \partial^2_{s_1} f_{\varepsilon}(u_1^o,u_2^o,t^o) \partial^2_{s_2} f_{\varepsilon}(u_1^o,u_2^o,t^o) - &[\partial_{s_1}\partial_{s_2} f_{\varepsilon}(u_1^o,u_2^o,t^o)]^2 \geq 0. 
\end{align}
The Young inequality, equation \ref{eq:rovna_se_2} and inequality \ref{eq:determinant} yield
\begin{align*}
\Delta f_{\varepsilon}(u_1^o,u_2^o,t^o) &= \partial^2_{s_1} f_\varepsilon(u_1^o,u_2^o,t^o) + \partial^2_{s_2} f_\varepsilon(u_1^o,u_2^o,t^o) \\
&\geq 2 [\partial^2_{s_1} f_\varepsilon(u_1^o,u_2^o,t^o) \; \partial^2_{s_2} f_\varepsilon(u_1^o,u_2^o,t^o)]^{\frac{1}{2}} \\
&\geq 2|\partial_{s_1} \partial_{s_2} f_\varepsilon(u_1^o,u_2^o,t^o)| = 4. 
\end{align*}
Using \ref{lem:heat} and the definition of $f_{\varepsilon}$ in \ref{eq:f_varepsilon} we obtain
\begin{align*}
0 &\geq \partial_t f_{\varepsilon}(u_1^o,u_2^o,t^o) = \partial_t f(u_1^o,u_2^o,t^o) + \varepsilon =\Delta f(u_1^o,u_2^o,t^o) - 4 + \varepsilon \\
&= \Delta f_{\varepsilon}(u_1^o,u_2^o,t^o) - 4 + \varepsilon \geq \varepsilon,
\end{align*}
which contradicts $\varepsilon > 0 $.
\end{proof}
\end{proposition}

The strength of \ref{thm:spherical} lies in the fact that the motion of each point on the curve is dictated by local information only, i.e. the normal vector and curvature, yet even parts of the curve separated by long distance along the curve are guarantied to avoid each other. Because of this local nature of the proof, the statement holds true even if the curve is separated into several disjoint closed curves which all simultaneously satisfy the assumptions.

Below, we introduce the notion of mutually spherical curves, formalizing this property.

\begin{definition}[Mutually spherical curves]
Family of space curves $\{\Gamma_{\alpha}\}_{\alpha \in \mathcal{M}}$ is mutually spherical if
\begin{align}\label{eq:def_mutually_spherical}
\exists x \in \mathbb{R}^3, \, \exists r > 0, \, \forall \alpha \in \mathcal{M}, \, \forall y \in \Gamma_{\alpha}: \; \Vert y - x \Vert = r. 
\end{align} 
\end{definition}

\begin{corollary}\label{cor:mutual}
Let each curve from family $\{\Gamma_{t,\alpha} \}_{t \in [0,t_{max}), \alpha \in \mathcal{M}}$ evolve according to the curve shortening flow given by \ref{eq:curvature_flow_1}\ref{eq:curvature_flow_2}. Assume that the initial condition $\{\Gamma_{0,\alpha} \}_{\mathcal{\alpha \in M}}$ is a mutually spherical and mutually disjoint family of simple curves and their curvatures $\kappa_\alpha(u,t)$ can be uniformly bounded by $C > 0$ for all $\alpha \in \mathcal{M}$, $u \in \mathbb{S}^1$ and $t \in [0,t_{max})$. Then the curves $\Gamma_{t,\alpha}$ cannot intersect for all $t \in (0,t_{max})$.
\end{corollary}

In \ref{cor:mutual}, we assume that all curves lie on the same sphere. However, similar result can be obtained for curves on different, mutually disjoint spheres. As the curves evolve, they remain embedded on spheres shrinking according to the mean curvature flow. Since the spheres remain disjoint, the curves cannot intersect either. We refer the reader to \cite{MinKimBen18} for further details.


\section{Conclusions}
\label{sec:conclusions}

This article contributes to the understanding of the long-term behavior of space curves during the curve shortening flow. We explored the differences between the flow in $\mathbb{R}^2$ and $\mathbb{R}^3$ and presented new properties for the latter. In \ref{thm:projectable}, we showed that even though space curves may stop being convex, the convexity of their orthogonal projection is preserved. Second part of the article is focused on the evolution of spherical curves and includes proof of the Avoidance principle (see \ref{thm:spherical}).

In practice, the theoretical results obtained in this work can help reduce computational time of numerical simulations. Knowing that the curves cannot intersect from the initial condition, we can switch off algorithmic treatment of topological changes. Note that such algorithms are highly time consuming as their computational complexity is usually $\mathcal{O}(n^2)$, where $n$ is the number of nodes on the discretized curve. For further reading see \cite{Balaz2012}, where an approach for reduction of the time complexity has been proposed for motion of curves in $\mathbb{R}^2$.

It remains an open question whether the technique of Gage and Hamilton can be used for proving the Avoidance principle for other families of space curves. Other promising research direction can be the study of knotted or linked parametric curves evolving in $\mathbb{R}^3$. The discussion can also be enriched by considering the problem in manifolds endowed with a Finsler metric or by increasing the dimension and/or codimension of the object in motion.


\printbibliography

\end{document}